\documentclass{amsart}
\usepackage[]{amsmath, amsthm, amsfonts, verbatim, amssymb}

\numberwithin{equation}{section}

\newcommand\C{{\mathbb C}}

\newcommand\dee{\partial}
\renewcommand\O{\Omega}
\newcommand\Obar{\overline{\Omega}}

\newcommand\Ohat{\widehat\Omega}

\begin{document}
 
\title[M\"obius transformations]
{M\"obius transformations, the Carath\'eodory metric, and the
objects of complex analysis and potential theory in multiply
connected domains}
\author[Steven R. Bell]
{Steven R. Bell}

\address[]{Mathematics Department, Purdue University, West Lafayette, IN  47907
USA} \email{bell@math.purdue.edu}

\thanks{Research supported by NSF grant DMS-0072197}



\begin{abstract}
It is proved that the family of Ahlfors extremal mappings
of a multiply connected region in the plane onto the unit disc
can be expressed as a rational combination of two fixed Ahlfors
mappings in much the same way that the family of Riemann mappings
associated to a simply connected region can be expressed in
terms of a single such map.  The formulas reveal that this
family of mappings extends to the double as a real analytic
function of both variables.  In particular, the infinitesimal
Carath\'eodory metric will be expressed in strikingly simple terms.
Similar results are proved for the Green's function, the Poisson
kernel, and the Bergman kernel.
\end{abstract}
 
\maketitle
 
\theoremstyle{plain}

\newtheorem {thm}{Theorem}[section]
\newtheorem {lem}{Lemma}

\hyphenation{bi-hol-o-mor-phic}
\hyphenation{hol-o-mor-phic}
\hyphenation{Car-a-the-o-dor-y}

\section{Introduction}
\label{intro}
Let $f_b$ denote the Riemann mapping function associated to
a point $b$ in a simply connected planar domain $\O\ne \C$.
Everyone knows that $f_b$ is the solution to an extremal
problem; it is the holomorphic map $h$ of $\O$ into the unit
disc such that $h'(b)$ is real and as large as possible.  Everyone
also knows that all the maps $f_b$ can be expressed in terms
of a single Riemann map $f_a$ associated to a point $a\in\O$
via
$$f_b(z)=\lambda\frac{f_a(z)-f_a(b)}{1-f_a(z)\overline{f_a(b)}}$$
where the unimodular constant $\lambda$ is given by
$$\lambda=\frac{\overline{f_a'(b)}}{|f_a'(b)|}.$$
In this paper, I shall prove that the solutions to the analogous
extremal problems on a finitely multiply connected domain in the plane,
the Ahlfors mappings, can be expressed in terms of just {\it two\/}
fixed Ahlfors mappings.  Many similarities with the
formula above in the simply connected case will become apparent and
I will explore some of the algebraic objects that present themselves.
A byproduct of these considerations will be that the infinitesimal
Carath\'eodory metric on a multiply connected domain is simply a
rational combination of two Ahlfors maps times one of their
derivatives.  I will explain an outlook which reveals a natural
way to view the extremal functions involved in the definition of the
Carath\'eodory metric ``off the diagonal'' in such a way that they
extend to $\Ohat\times\Ohat$ where $\Ohat$ is the double of $\O$.

I will also investigate the complexity of the classical Green's
function and Bergman kernel associated to a multiply connected domain.
In particular, it is proved in \S\ref{sec6} that if $\O$ is a
finitely connected domain in the plane such that no boundary component
is a point, then  there exist two Ahlfors maps $f_a$ and $f_b$ associated
to $\O$ such that the Bergman kernel for $\O$ is given by
$$K(w,z)=
\frac{f_a'(w)\overline{f_a'(z)}}{(1-f_a(w)\overline{f_a(z)})^2}
\left(\sum_{j,k=1}^N H_j(w)\overline{H_k(z)}\right),$$
where the functions $H_j$ are rational combinations of the
two Ahlfors maps $f_a$ and $f_b$.  Future avenues of research include
the problem of extending these results to finite Riemann surfaces and
the problem to determine the way the rational functions that arise in
these formulas depend on the domain.

\section{The smooth case}
\label{sec2}
To get started, we shall assume that $\O$ is a bounded
$n$-connected domain in the plane with $C^\infty$ smooth boundary consisting
of $n$~non-intersecting curves.  (Later, we shall consider general
$n$-connected domains such that no boundary component is a point.)
Let $S(z,w)$ denote the Szeg\H o kernel associated to $\O$ (see
\cite{r2} or \cite{r8} for definitions and standard terminology in
what follows).

Fix a point $a$ in $\O$ so that the $n-1$~zeroes $a_1,\dots,a_{n-1}$ of
$S(z,a)$ in the $z$-variable are distinct simple zeroes.  (That such
points $a$ form an open dense subset of $\O$ was proved in \cite{r3}.)
Let $a_0$ be equal to $a$.  I proved in \cite[Theorem~3.1]{r4}
that the Szeg\H o kernel can be expressed in terms of the $n+1$
functions of one variable, $S(z,a)$, $f_a(z)$, and $S(z,a_i)$, $i=1,\dots,n-1$
via the formula
\begin{equation}
S(z,w)=\frac{1}{1-f_a(z)\overline{f_a(w)}}
\sum_{i,j=0}^{n-1} c_{ij}S(z,a_i)\,\overline{S(w,a_j)}
\label{e2.1}
\end{equation}
where $f_a(z)$ denotes the Ahlfors map associated to $(\O,a)$
and the coefficients $c_{ij}$ are given as the
coefficients of the inverse matrix to the $n\times n$ matrix
$\left[S(a_j,a_k)\right]$.
A similar formula for the Garabedian kernel was proved in \cite{r5},
\begin{equation}
L(z,w)=\frac{f_a(w)}{f_a(z)-f_a(w)}
\sum_{i,j=0}^{n-1} c_{ij}S(z,a_i)L(w,a_j),
\label{e2.2}
\end{equation}
where the constants $c_{ij}$ are the same as the constants in (\ref{e2.1}).

Given a point $w\in\O$, the Ahlfors map $f_w$ associated to the pair $(\O,w)$
is a proper holomorphic mapping of $\O$ onto the unit disc.  It is an
$n$-to-one mapping (counting multiplicities), it extends to be in
$C^\infty(\Obar)$, and it maps each boundary curve of $\O$ one-to-one
onto the unit circle.  Furthermore, $f_w(w)=0$, and $f_w$ is the unique
function mapping $\O$ into the unit disc maximizing the quantity $|f_w'(w)|$
with $f_w'(w)>0$.  The Ahlfors map is related to the Szeg\H o kernel
and Garabedian kernel via (see \cite[page~49]{r2})
\begin{equation}
f_w(z)=\frac{S(z,w)}{L(z,w)}.
\label{e2.3}
\end{equation}

When equations~(\ref{e2.1}) and (\ref{e2.2}) are substituted into
(\ref{e2.3}), we obtain the monstrosity

\begin{equation*}
f_w(z)=
\frac{f_a(z)-f_a(w)}
{f_a(w)(1-f_a(z)\overline{f_a(w)}\,)}
\ \frac{
\sum_{i,j=0}^{n-1} c_{ij}S(z,a_i)\,\overline{S(w,a_j)}
}{
\sum_{i,j=0}^{n-1} c_{ij}S(z,a_i)L(w,a_j)
}.
\end{equation*}

Next, divide the numerator and the denominator of the second quotient
in this expression by $S(z,a)S(w,a)$ and multiply the whole thing by
one in the form of $\overline{S(w,a)}/\overline{S(w,a)}$ to obtain

\begin{equation*}
f_w(z)=
\frac{f_a(z)-f_a(w)}
{f_a(w)(1-f_a(z)\overline{f_a(w)}\,)}
\left(
\frac{
\sum_{i,j=0}^{n-1} c_{ij}\frac{S(z,a_i)}{S(z,a)}\,
{\overline{S(w,a_j)}}/{\overline{S(w,a)}}
}{
\sum_{i,j=0}^{n-1} c_{ij}\frac{S(z,a_i)}{S(z,a)}\,
{L(w,a_j)}/{S(w,a)}
}
\right)
\frac{\overline{S(w,a)}}{S(w,a)}.
\end{equation*}

It is not hard to show that $f_a(z)$ and quotients of the form
$S(z,a_i)/S(z,a)$ and $L(z,a_i)/S(z,a)$ extend to the double $\Ohat$
of $\O$ as meromorphic functions (see \cite[p. 6]{r6}).  Since the
argument is quick and simple, we give it here.  Let $R(z)$ denote the
antiholomorphic reflection function on $\Ohat$ which maps $\O$ into
the reflected copy of $\O$.  Note that $f_a(z)$ is equal to
$1/\overline{f_a(z)}$ on $b\O$, which is equal to
$1/\overline{f_a(R(z))}$ there.  Hence, the holomorphic
function $f_a(z)$ on $\O$ and the meromorphic function
$1/\overline{f_a(R(z))}$ on the complement of $\O$ in $\Ohat$ both
extend continuously up to $b\O$ and have the same values there.  Hence
$f_a$ extends meromorphically to the double.  Similar reasoning can
be applied to the quotients as follows.  The Garabedian kernel is
related to the Szeg\H o kernel via the identity
\begin{equation}
\frac{1}{i} L(z,a)T(z)=S(a,z)\qquad\text{for $z\in b\O$ and $a\in\O$,}
\label{e2.4}
\end{equation}
where $T(z)$ denotes the complex number of unit modulus pointing in
the tangent direction at $z\in b\O$ chosen so that $iT(z)$ represents
an inward pointing normal vector to the boundary.  Hence,
$S(z,a_i)/S(z,a)$ is equal to the conjugate of $L(z,a_i)/L(z,a)$ on
the boundary and the same reasoning used above for $f_a$ shows that
$S(z,a_i)/S(z,a)$ extends to the double meromorphically.  Similarly
$L(z,a_i)/S(z,a)$ is equal to the conjugate of $S(z,a_i)/L(z,a)$ on
the boundary and this shows that $L(z,a_i)/S(z,a)$ extends to the
double meromorphically.

It is proved in \cite{r7} that it is possible to choose a second
Ahlfors map $f_b$ so that $f_a$ and $f_b$ generate the field of
meromorphic functions on $\Ohat$.  (Such a pair is called a {\it
primitive pair}, see \cite{r1} and \cite{r9}).  
Hence, we have now shown that there exists a rational function on
$\C^6$ such that
\begin{equation}
f_w(z)=
\lambda(w)R(f_a(z),f_b(z),f_a(w),f_b(w),\overline{f_a(w)},\overline{f_b(w)})
\label{e2.5}
\end{equation}
where $\lambda(w)$ is the unimodular function given by
$$\lambda(w)=\overline{S(w,a)}/S(w,a).$$
This formula is very reminiscent of the formula for the Riemann
maps mentioned at the beginning of this paper.
It now becomes irresistible to drop the factor $\lambda(w)$ from
equation~(\ref{e2.5}) and to define a function $F(z,w)$ via
$$F(z,w)=f_w(z)/\lambda(w).$$
Let us call this function the {\it alternatively normalized Ahlfors map}.
Under this normalization, the map $z\mapsto F(z,w)$ has a derivative
at $w$ with extremal modulus, however the argument of the derivative
is $-\arg \lambda(w)$ there instead of zero.  This family of extremal
maps has the astonishing feature that it extends in a unique way to
$\Ohat\times\Ohat$ as a complex rational function of $f_a(z)$, $f_b(z)$, and
$f_a(w)$, $f_b(w)$, $\overline{f_a(w)}$, $\overline{f_b(w)}$.  Furthermore,
this extension is meromorphic in $z$ and real analytic in $w$.  One
might also glimpse some semblance of an analogue of a M\"obius function
in these deliberations and we shall come back to this point later in
the paper.

Another important consequence of formula (\ref{e2.5}) is that the
infinitesimal Carath\'eo\-dory metric can be expressed in terms of two
Ahlfors maps.  In fact, it is shown in \cite[p.~344]{r7} that the quotient
$f_b'(z)/f_a'(z)$ extends to be meromorphic on the double of $\O$
and is therefore a rational combination of $f_a(z)$ and $f_b(z)$.
Hence, if we differentiate~(\ref{e2.5}) with respect to $z$ and take
the modulus of the expression, we obtain that $|f_w'(z)|$ is given
by $|f_a'(z)|$ times the modulus of a rational function of
$f_a(z)$, $f_b(z)$, $f_a(w)$, $f_b(w)$, $\overline{f_a(w)}$, and
$\overline{f_b(w)}$.  Now, if we set $w=z$, we may conclude that the
infinitesimal Carath\'eodory
metric is given by $\rho(z)|dz|$ where
$$\rho(z)=|f_a'(z)|
\left|
Q(f_a(z),f_b(z),\overline{f_a(z)},\overline{f_b(z)})
\right|,$$
where $Q$ is a rational function on $\C^4$.

Many questions present themselves at this point.  The formula above
for the infinitesimal Carath\'eodory metric almost looks exact.  Might
there exist special multiply connected domains where the
Carath\'eodory metric could be computed as easily as it is in the
unit disk?  Another natural question to ask is whether or not similar
formulas hold for finite Riemann surfaces.  Ahlfors mappings are
available in this setting, but the relationship between these maps and
the kernel functions used in the proof in the planar case are not as
straightforward.  New methods of proof would have to be discovered.

\section{The nonsmooth case}
\label{sec3}
Suppose that $\O$ is merely an $n$-connected domain in the plane such
that no boundary component is a point.  It is well known that there is
a biholomorphic mapping $\phi$ mapping $\O$ one-to-one onto a bounded domain
$\O_a$ in the plane with smooth real analytic boundary.  The standard
construction yields a domain $\O_a$ that is a bounded $n$-connected domain with
$C^\infty$ smooth boundary
whose boundary consists of $n$ non-intersecting simple closed real analytic
curves.  Let subscript or superscript $a$'s indicate that a kernel function
or mapping is associated to $\O_a$.  Kernels without sub or superscripts are
associated to $\O$.  It is well known that the function $\phi'$ has a single
valued holomorphic square root on $\O$ (see \cite[p. 43]{r2}).  We
define the Szeg\H o kernel and Garabedian kernel associated to $\O$ via
the natural transformation formulas,
\begin{equation*}
S(z,w)=\sqrt{\phi'(z)}\ S_a(\phi(z),\phi(w))\overline{\sqrt{\phi'(w)}}
\end{equation*}
and
\begin{equation*}
L(z,w)=\sqrt{\phi'(z)}\ L_a(\phi(z),\phi(w))\sqrt{\phi'(w)}.
\end{equation*}
The Ahlfors map associated to a point $b\in\O$ is defined to be
the solution to the extremal problem, $f_b:\O\to D_1(0)$ with $f_b'(b)>0$ and
maximal.  It is easy to see that Ahlfors maps satisfy
\begin{equation*}
f_b(z)=\lambda f_{\phi(b)}^a(\phi(z))
\end{equation*}
for some unimodular constant $\lambda$ and it follows that
$f_b(z)$ is a proper holomorphic mapping of $\O$ onto $D_1(0)$.
It also follows that $f_b(z)$ is given by $S(z,b)/L(z,b)$ just as in
the smooth case.  Now it is easy to see that all the quotients that
appeared in the proofs of the results of \S\ref{sec2} are invariant
under $\phi$ and the proofs carry over line for line.  We may now
state the following theorem.

\begin{thm}
Suppose that $\O$ is an $n$-connected domain in the plane such
that no boundary component is a point.  There exist two points $a$ and
$b$ in $\O$ such that the alternatively normalized Ahlfors map $F(z,w)$
associated to $\O$ is a complex rational function of $f_a(z)$,
$f_b(z)$, and $f_a(w)$, $f_b(w)$, $\overline{f_a(w)}$, $\overline{f_b(w)}$.
Furthermore, the family of Ahlfors mappings is given by formula~(\ref{e2.5})
and the infinitesimal Carath\'eodory metric is given by $\rho(z)|dz|$ where
$$\rho(z)=|f_a'(z)|
\left|
Q(f_a(z),f_b(z),\overline{f_a(z)},\overline{f_b(z)})
\right|,$$
where $Q$ is a rational function on $\C^4$.
\end{thm}

\section{What is a M\"obius transformation?}
\label{sec4}
Here is one way to ``invent'' M\"obius transformations.  Let $p(z)$
denote an irreducible polynomial of one variable with no zeroes
in the unit disc, i.e., let $p(z)=z-b$ where $|b|>1$.  Notice
that $p(1/\bar z)$ is equal to $p(z)$ on the unit circle.  Let
$q(z)$ denote the polynomial obtained by multiplying the conjugate
of $p(1/\bar z)$ by the power of $z$ needed to clear the poles
in the unit disc, i.e., $q(z)=1-z\bar b$.  Since
$|q(z)|=|\bar zp(1/\bar z)|$, it follows that $|q(z)|=|p(z)|$
on the unit circle.  Notice that $q(z)/p(z)$ is a M\"obius transformation
(let $b=1/\bar a$ to make it look more standard).

It is shown in \cite{r7} that every proper holomorphic mapping of
a smooth $n$-connected domain $\O$ onto the unit disk can be expressed
as a rational combination of two Ahlfors maps $f_a$ and $f_b$
associated to points $a$ and $b$ in $\O$.  It is an interesting
problem to determine just exactly which rational functions arise
in this manner, and it is tempting to call some of these rational
functions M\"obius transformations.  Here is one way to construct
such a rational function.  Let $\Delta^2$ denote the unit bidisc.  Let
$p(z,w)$ denote an irreducible polynomial of two variables with no
zeroes in the closure of $\Delta^2$.  Notice
that $p(1/\bar z,1/\bar w)$ is equal to $p(z,w)$ on the distinguished
boundary of $\Delta^2$.  Suppose $N$ is the degree of $p(z,w)$ in $z$
and $M$ is the degree in $w$.  Let $q(z,w)$ be the polynomial given
by $z^Nw^M$ times the conjugate of $p(1/\bar z,1/\bar w)$.  Since
$q(z,w)$ and $p(z,w)$ have the same modulus on the distinguished boundary
of $\Delta^2$, and since $|z^Nw^M|=1$ there, it follows that the modulus of
$q(z,w)/p(z,w)$ is also one there.   It follows that, if it is not constant,
$q(f_a(z),f_b(w))/p(f_a(z),f_b(w))$ is a proper holomorphic mapping
of $\O$ onto the unit disc.

More generally, the same construction can be carried out if $p(z,w)$ is an
irreducible polynomial on $\C^2$ which does not vanish on the portion
of the curve $z\mapsto (f_a(z),f_b(z))$ inside the closed unit
bidisc.  Can any proper map from $\O$ to the unit disc be expressed in a
similar manner, perhaps as some kind of combination of these basic
maps?

\section{The Poisson kernel extends nicely to the double}
\label{sec5}
Of course the Poisson kernel extends to the double by simple
reflection.  Here we show that it extends nicely in both variables and
in terms of some special functions with geometric meaning.

Assume that $\O$ is a bounded $n$-connected domain in the plane with
$C^\infty$ smooth boundary consisting of $n$~non-intersecting curves.
Let $\gamma_1,\dots,\gamma_{n-1}$ denote the inner curves and let
$\gamma_n$ denote the outer curve.

The classical Poisson kernel for $\O$ is related to the normal derivative
of the Green's function via
$$p(z,w)=\frac{1}{2\pi}\frac{\dee}{\dee n_w}G(z,w)\qquad z\in\O,\ w\in
b\O,$$
where $(\dee/\dee n_w)$ denotes the normal derivative in the $w$
variable.  It is a standard fact that we may rewrite this last formula
(see \cite[pages~134-136]{r2}) in the form
$$p(z,w)=-\frac{i}{\pi}\frac{\dee}{\dee w}G(z,w){T(w)}.$$
It is proved in \cite[p. 1367]{r4} (see also \cite[p. 12]{r6} for an
easier proof) that the derivative of the Green's function
$G_w(z,w):=\frac{\dee}{\dee w}G(z,w)$ is given by
\begin{equation}
G_w(z,w)=
\pi\frac{S(w,z)L(w,z)}{S(z,z)}
+i\pi\sum_{j=1}^{n-1}
\left(\omega_j(z)-\lambda_j(z)\right)u_j(w),
\label{e5.1}
\end{equation}
where the functions $\lambda_j(z)$ are given by
$$\lambda_j(z)=
\int_{w\in\gamma_j}\frac{|S(w,z)|^2}{S(z,z)}\ ds,$$
the functions $\omega_j(z)$ are the harmonic measure
functions, and the functions $u_j$ are a basis for the
linear span $\mathcal F'$ of the functions $F_j':= 2(\dee\omega_j/\dee z)$
normalized so that
$$\delta_{kj}=\int_{\gamma_k}u_j(w)\ dw.$$

We now show that the principal term
$\frac{S(w,z)L(w,z)}{S(z,z)}$ in the expression for $G_w(z,w)$
has the interesting property that it extends to the double of
$\O$ in the $z$ variable as a real analytic function that is a
rational combination of two Ahlfors maps $f_a(z)$ and $f_b(z)$
and their conjugates.  Indeed, if we substitute equations~(\ref{e2.1})
and (\ref{e2.2}) into this expression, we obtain that
\begin{equation}
\frac{S(w,z)L(w,z)}{S(z,z)}= T_1(z,w) T_2(z,w)
\label{e5.2}
\end{equation}
where
$$
T_1(z,w)=
\frac{(1-|f_a(z)|^2)f_a(z)}
{(f_a(w)-f_a(z))(1-f_a(w)\overline{f_a(z)}\,)}
$$
and
$$
T_2(z,w)=
\frac{\left(
\sum_{i,j=0}^{n-1} c_{ij}S(w,a_i)\,\overline{S(z,a_j)}
\right)
\left(
\sum_{i,j=0}^{n-1} c_{ij}S(w,a_i)L(z,a_j),
\right)}
{\sum_{i,j=0}^{n-1} c_{ij}S(z,a_i)\,\overline{S(z,a_j)}}.
$$
The first term extends to the double as a real analytic function
because $f_a$ does.  If we divide the numerator and denominator of
the second term by $|S(z,a)|^2$, we observe that the numerator is
a linear combination of functions which are given as products of
$L(z,a_n)/S(z,a)$ times
the conjugate of
$S(z,a_m)/S(z,a)$ times
$S(w,a_q)S(w,a_k)$.
As mentioned previously, the functions
$S(z,a_m)/S(z,a)$ and $L(z,a_n)/S(z,a)$ extend meromorphically to the
double, and hence can be expressed as rational functions of two Ahlfors
maps $f_a(z)$ and $f_b(z)$.
The denominator is a linear combination of functions given as the
product of
$L(z,a_n)/S(z,a)$ times
the conjugate of $S(z,a_m)/S(z,a)$.
Hence, it has these properties, too, in the $z$ variable.

We have shown that, for fixed $w$, the function of $z$ given by
$S(w,z)L(w,z)/S(z,z)$ is a rational combination of two
Ahlfors maps $f_a(z)$ and $f_b(z)$ and their conjugates.

We now claim that functions of the form
$$S(w,a_q)S(w,a_k)/f_a'(w)$$
extend meromorphically to the double of $\O$.  Indeed, since
identity~(\ref{e2.4}) yields that
$S(w,a_q)S(w,a_k)T(w)$ is equal to the conjugate of
$-L(w,a_q)L(w,a_k)T(w)$ for $w$ in the boundary and since
$T(w)f_a'(w)/f_a(w)$ is equal to the conjugate of
$-T(w)f_a'(w)/f_a(w)$, we may use similar reasoning to that in
\cite[p. 202]{r5} and divide these two expressions to see that
$S(w,a_q)S(w,a_k)f_a(w)/f_a'(w)$ extends meromorphically to the
double.  Since $f_a(w)$ extends to the double, the claim is proved.
Now, if we were to divide the large expression for
$S(w,z)L(w,z)/S(z,z)$ above by $f_a'(w)$, we would
deduce that
$$\frac{S(w,z)L(w,z)}{f_a'(w)S(z,z)}$$
is a rational combination of the two functions $f_a(w)$ and $f_b(w)$, and
the four functions
$f_a(z)$ and $f_b(z)$ and their conjugates.

It is proved in \cite{r6} that there exist
$n-1$ points $w_j$ in $\O$ such that the functions
$$(\omega_k(z)-\lambda_k(z))$$
are linear combinations of functions of $z$ of the form
$$G_w(z,w_j)-\pi\frac{S(w_j,z)L(w_j,z)}{S(z,z)}.$$
The function $G_w(z,w_j)$ is harmonic in $z$ on $\O-\{w_j\}$
and vanishes on the boundary.  Hence, it extends to the double
as a harmonic function with two singular points.
The function
$S(w_j,z)L(w_j,z)/S(z,z)$
has been shown to extend to the double.
Since $F_j'(w)T(w)$ is equal to the conjugate of
$-F_j'(w)T(w)$ on the boundary, the same reasoning used above yields
that functions of the form
$$F_j'(w)/f_a'(w)$$
extend meromorphically to the double of $\O$.
We may now state that $G_w(z,w)$ is given as $f_a'(w)$ times a
rational combination of the two functions $f_a(w)$ and $f_b(w)$, and
the four functions $f_a(z)$ and $f_b(z)$ and their conjugates plus
a linear combination of the functions $G_w(z,w_j)$ times rational
combinations of $f_a(w)$ and $f_b(w)$.  In symbols,
\begin{eqnarray*}
G_w(z,w) &=& f_a'(w)
R_0(f_a(w),f_b(w),f_a(z),f_b(z),\overline{f_a(z)},\overline{f_b(z)}\,)
\nonumber\\
& &{}+f_a'(w)\sum_{j=1}^{n-1} G_w(z,w_j)R_j(f_a(w),f_b(w))\nonumber
\end{eqnarray*}
where the functions $R_0$ and $R_j$ are rational.  All the functions
that comprise $G_w(z,w)$ extend nicely to the double except
$f_a'(w)$.

The results of this section can be generalized to $n$-connected domains
with non-smooth boundaries in the same way as in \S\ref{sec3}, but
we shall not do this here.

\section{Linearizing the Green's function and Bergman kernel}
\label{sec6}
In the simply connected case, the Green's function is related to
a Riemann map $f(z)$ by the very simple formula
$$G(z,w)=\ln\left|\frac{f(z)-f(w)}{1-f(z)\overline{f(w)}}\right|.$$
In the multiply connected setting, the Green's function is also related
to Ahlfors maps, but it is not clear if the Green's function can be
expressed naturally in terms of maps.  We saw some tantalizing
evidence in the last section that there might be such an expression.
In this section, I give some further evidence which leads me to believe
that maybe such an expression exits.  This evidence fits nicely into the
subject matter of this paper because a genuine M\"obius transformation
is a key ingredient.

Suppose that $\O$ is a multiply connected domain with $C^\infty$-smooth
boundary and let $f(z)$ denote an Ahlfors map associated to $(\O,a)$
that has simple zeroes.  Let $\mathcal L(z,w)$ denote the function
$$\ln\left|\frac{f(z)-f(w)}{1-\overline{f(z)}f(w)}\right|.$$
We want to investigate the boundary behavior of the quotient
$G(z,w)/\mathcal L(z,w)$ as $z$ and $w$ are both allowed to approach
the boundary.  First assume that $z$ is a fixed point in $\O$ and let $w$
approach the boundary.  Both the numerator and the denominator extend
$C^\infty$ smoothly in $w$ to the boundary and the Hopf Lemma reveals
that both terms vanish to first order along the boundary.  Hence, the quotient
extends $C^\infty$-smoothly up to the boundary in $w$ and the limit is
given (by L'H\^opital's rule) as the quotient of the normal derivatives
$(\dee/\dee n_w)$ in the $w$ variable.  Recall that
$$\frac{\dee}{\dee n_w}G(z,w)= -2iG_w(z,w)T(w).$$
Since $\mathcal L$ is also a real valued harmonic function that vanishes
on the boundary, the same reasoning that yields this identity can be
applied to the normal derivative of $\mathcal L(z,w)$ to obtain
$$\frac{\dee}{\dee n_w}\mathcal L(z,w)= -2i\mathcal L_w(z,w)T(w)=
\frac{f'(w)(1-|f(z)|^2)T(w)}{i(f(w)-f(z))(1-\overline{f(z)}f(w))}.$$
Notice the similarity of this expression with $T_1(z,w)$ in
formula~(\ref{e5.2}).  We may now divide these two normal derivatives and use
(\ref{e5.1}) and (\ref{e5.2}) to obtain that
\begin{equation*}
\frac{(\dee/\dee n_w)G(z,w)}
{(\dee/\dee n_w)\mathcal L(z,w)}=
\frac{if(z)}{f'(w)}T_2(z,w) + T_3(z,w)
\end{equation*}
where
$$
T_2(z,w)=
\frac{\left(
\sum_{i,j=0}^{n-1} c_{ij}S(w,a_i)\,\overline{S(z,a_j)}
\right)
\left(
\sum_{i,j=0}^{n-1} c_{ij}S(w,a_i)L(z,a_j),
\right)}
{\sum_{i,j=0}^{n-1} c_{ij}S(z,a_i)\,\overline{S(z,a_j)}}.
$$
and
$$
T_3(z,w)=
\frac
{i(f(w)-f(z))(1-\overline{f(z)}f(w))}
{f'(w)(1-|f(z)|^2)}
\sum_{j=1}^{n-1} i\pi
\left(\omega_j(z)-\lambda_j(z)\right)u_j(w).$$
Although this formula is painful to look at, a moment
of suffering reveals that the right hand side can
be written as a sum of simple terms to yield that
\begin{equation}
\frac{(\dee/\dee n_w)G(z,w)}
{(\dee/\dee n_w)\mathcal L(z,w)}=
\sum_{j}\mu_j(z)h_j(w)
\label{e6.1}
\end{equation}
where each of the functions $\mu_j$ extend to the double as
real analytic functions and each of the functions $h_j(w)$ extend
to the double as meromorphic functions.  (Note that here we have
used the fact proved earlier that $u_j/f'$ extends to the double
as a meromorphic function.)  I view this formula as a linearization
or polarization of the Poisson kernel.  I take this opportunity to
state a theorem.

\begin{thm}
Suppose that $\O$ is a bounded finitely connected domain in the plane
with $C^\infty$ smooth boundary.  The Green's function associated to
$\O$ satisfies an identity of the form
\begin{equation*}
\frac{\dee}{\dee w}G(z,w)= \frac{\dee}{\dee w}\mathcal L(z,w)
\sum_{j}\mu_j(z)h_j(w)
\end{equation*}
where each $\mu_j$ extends to be real analytic on the double of $\O$
and each $h_j$ extends to be meromorphic on the double of $\O$.
\label{thm6.1}
\end{thm}

We note that we have proved this identity when $z$ is in $\O$ and $w$
is in the boundary of $\O$, but since the functions of $w$ in the expression
are all meromorphic, the identity extends to hold for all $w$ in $\Obar$.

We now continue to deal with equation~(\ref{e6.1}) and we assume that
$w$ is back in the boundary and we let the $z$ variable tend to a
boundary point other than $w$ to obtain
\begin{equation*}
\frac{(\dee^2/\dee n_z\dee n_w)G(z,w)}
{(\dee^2/\dee n_z\dee n_w)\mathcal L(z,w)}=
\frac{(\dee^2/\dee \bar z\dee w)G(z,w)}
{(\dee^2/\dee \bar z\dee w)\mathcal L(z,w)}=
\frac{if(z)}{f'(w)}T_2(z,w) + T_4(z,w),
\end{equation*}
where $T_2(z,w)$ is as above, but $T_4(z,w)$ is given by
$$
T_4(z,w)=
\frac
{i(f(w)-f(z))(1-\overline{f(z)}f(w))}
{f'(w)}\sum_{j=1}^{n-1}\nu_j(z)u_j(w),$$
where $\nu_j(z)$ is equal to the limit of
$i\pi\left(\omega_j(z)-\lambda_j(z)\right)/(1-|f(z)|^2)$
as $z$ tends to the boundary.
Since $(\dee^2/\dee \bar z\dee w)G(z,w)$ is equal to $K(w,z)$ and since
$(\dee^2/\dee \bar z\dee w)\mathcal L(z,w)$ is equal to
$$\frac{f'(w)\overline{f'(z)}}{(1-f(w)\overline{f(z)})^2},$$
we deduce that
$$K(w,z)=
\frac{f'(w)\overline{f'(z)}}{(1-f(w)\overline{f(z)})^2}
\left(\sum_j \sigma_j(z)H_j(w)\right)
$$
where the sum is finite and each function $H_j(w)$ extends to be
meromorphic on the double of $\O$.  We may assume that this sum has
been collapsed so that the functions $H_j(w)$ are linearly independent
on $\O$.   We can now exploit the hermitian property of the Bergman
kernel to easily deduce that the $\sigma_j$ functions are actually
linear combinations of the conjugates of the $H_j$.  Hence, we have
proved that
\begin{equation*}
K(w,z)=
\frac{f'(w)\overline{f'(z)}}{(1-f(w)\overline{f(z)})^2}
\left(\sum_{j,k=1}^N H_j(w)\overline{H_k(z)}\right).
\end{equation*}
We have only shown that this identity holds on the boundary, but
it is clear that it extends to the inside of the domain because
all the functions that appear in the identity are meromorphic.
We can now use the fact proved in \cite{r7} that the field of
meromorphic functions on the double is generated by two Ahlfors
maps to be able to state that the functions $H_j$ are rational
combinations of two Ahlfors maps.  We have operated under the
assumption that $\O$ has smooth boundary.  Finally, if $\O$ does
not have smooth boundary, we can map to a domain with smooth
boundary and use the fact that the terms in the expression for
$K(z,w)$ transform under biholomorphic mappings to obtain
the following theorem.

\begin{thm}
Suppose that $\O$ is a finitely connected domain in the plane
such that no boundary component is a point.  There exist two
points $a$ and $b$ in $\O$ such that the Bergman kernel
associated to $\O$ is given by
$$K(w,z)=
\frac{f_a'(w)\overline{f_a'(z)}}{(1-f_a(w)\overline{f_a(z)})^2}
\left(\sum_{j,k=1}^N H_j(w)\overline{H_k(z)}\right),$$
where the functions $H_j$ are rational combinations of the
two Ahlfors maps $f_a$ and $f_b$.
\label{thm6.2}
\end{thm}

There are many interesting questions that present themselves
at this point.  The rational functions that appear in the formula
in Theorem~\ref{thm6.2} most likely satisfy an invariance property under
biholomorphic mappings and have algebraic geometric significance.
The functions $\lambda_j$ have many interesting properties.
I wonder if they might be expressible as rational combinations of
two Ahlfors maps and their conjugates.  I also wonder if the
Green's function can be shown to have similar finite complexity
to all the other kernel functions that have been studied in this
paper, modulo some logarithmic expressions.  It is also a safe
bet that many of the results in this paper extend to the case of
finite Riemann surfaces.  I leave these investigations for the future.

\end{document}